# Total Lagrangian Finite Element Formulation of the Flory-Rehner Free Energy Function

# Formulación de Elemento Finito Total Lagrangiano de la Función de Energía Libre de Flory-Rehner

Mario J. Juha[*]

Universidad Autónoma del Caribe

Barranquilla - Colombia

**Abstract**

We address the total Lagrangian finite element implementation of the Flory-Rehner free-energy function in the framework of a hyperelastic material model. We explicitly give all the equations required to implement this material model in an implicit nonlinear finite element analysis, particularly, we show how to derive the so-called algorithmic or consistent linearized tangent modulli in the Lagrangian description. Some analytical and numerical results for different boundary-value problems are presented to validate the implementation.

**Resumen**

Se trata la implemetación total Lagrangiana en elemento finito de la función de energía libre de Flory-Rehner en un marco de un modelo de material hiperelastico. Explícitamente se dan todas las ecuaciones requeridas para implementar este modelo de material en un análisis de elemento

---

[*] Assistant Research Professor, Department of Mechanical Engineering. Ph.D in Civil Engineering, Structures, University of South Florida (USF). B.Sc. Mechanical Engineering, Universidad del Norte. mario.juha@uac.edu.co. Tel: +57 5 3366800 Ext. 993

finito no lineal, particularmente, se muestra como derivar el llamado modulo tangente consistente o algorítmico en una descripción Lagrangiana. Algunos resultados analíticos y numéricos para diferentes problemas de valor de frontera son presentados para validar la implementación.

**Palabras clave:** Flory-Rehner, energía libre, total Lagrangiano, elemento finito, geles

# 1. INTRODUCTION

Flory-Rehner gels are a kind of soft matter that swell/shrink spontaneously due to a change of the chemical potential in the environment [1]. These gels are composed of a three dimensional polymeric network and small solvent molecules that can migrate in and out of the polymeric chains. The bonds between solvent and polymeric molecules can be thought as physical bonds on the order of Van der Waals forces [2]. Flory-Rehner gels are highly compressible elastic solids with a non linear material behavior that undergoes large deformations when they are in contact with a solvent at a different chemical potential. Therefore, due to these nonlinearities, it is not an easy task to develop simple theories to analyze them for general cases. Analysis of swelling of inhomogeous gels has been an active topic of study that dates back, at least, since the pioneering work of Tanaka and Fillmore [3] in the late 70's, when they tried to estimate, analytically, the characteristic time of swelling. Of particular interest is the mathematical work of Onuki [4, 5]. More recently, Dolbow et al. [6, 7] considered a continuum model for chemically induced volume phase transitions in hydrogels that allows for a sharp interface separating swelled and collapsed phases for the underlying polymer network. They used a hybrid eXtended-Finite Element/Level-Set method for obtaining approximate solutions to the governing equations for a two-dimensional model and get an estimate of the response time for swelling. Most notable is the recent work done by the group lead by Zhigang Suo [8, 9, 10], they presented a simple theory of coupled molecular migration and large deformation within the framework of non equilibrium thermodynamics where all the fields are given as a function of the material coordinates. They implemented a hyperelastic material model using the user supplied routine UHYPER in the commercial package ABAQUS. Kang and Huang [11, 12] have reformulated and implemented the previous theory using the user supplied routine UMAT of ABAQUS to analyze surface instabilities of confined gels layers. The intent of this paper is to address the continuum mechanics theory and finite element implementation of a Flory-Rehner gel in a total Lagrangian framework. Although the theory presented in this paper is not new, still it has not been fully studied and to the knowledge of the authors, there is not a reference where explicitly a linearization of the internal nodal forces have been done in a total Lagrangian formulation for a Flory-Rehner gel. However, in [11] a tangent modulus tensor for the same material studied here was explicitly presented for current configuration. Someone could argue that a series of pull-back operation [13] can be applied to the tangent modulus in the deformed representation such that we could transform it into the total Lagrangian formulation. Nevertheless, we have decided to present the entire structure of the total Lagrangian implementation. Following same lines as in [8, 14], we presented a purely mechanical theory of coupled species transport and large elastic deformation in a continuum mechanics approach.

## 2. GENERAL SETUP

Let introduce a fixed right-handed Cartesian coordinate system as presented in figure 1, with orthogonal unit base vectors, $E_i, e_i, i=1,2,3$ for the reference and current configuration, respectively. Select $\Omega_0$ as the volume occupy by the body in the reference state and $\Omega$ as the volume occupy by the body in the current state. We will assume that there exist a one-to-one and smooth mapping, $X=\chi^{-1}(x,t)$ , between the reference and current configuration for all time, $t \geq 0$ . The tensor field $F = \nabla \chi$ is referred as the deformation gradient and should satisfy the requirement that $J = det\, F > 0$ for all $t$ . The latter guarantee that the inverse motion, $X = \chi^{-1}(x,t)$ exist and it is unique. The boundary of the body is divided into three parts as follows: $\Gamma_h \cup \Gamma_t \cup \Gamma_u = \emptyset$ and with the restriction that $\Gamma = \Gamma_h \cup \Gamma_t \cup \Gamma_u$ . In the following subsections we will review the balance principles.

**Molecular incompressibility for gels**

The gel is a continuum body that could have a dramatic change in volume, therefore it is considered a compressible solid. If we look independently at the long chain network of the polymer and the solvent molecules that constitute the gel we will see that they are incompressible. Then, taking into account the incompressibility of the polymer and solvent molecules, and according to [14], the molecular incompressibility is defined in equation (1).

$$det(F) = 1 + vC \qquad (1)$$

Where $vC$ is the volume of the small molecules in the gel divided by the volume of the dry network. In figure 2 we present a pictorial definition of the hydrogel. Note that, although, the system is incompressible, the hydrogel is compressible, because the solvent particles could move in and out of the network. The kinetics of the swelling gels is presented in a pictorial form in figure 3. In the theory presented in this section, we have studied the state when the gel is in equilibrium with the solvent. In this steady state, there is no change in volume, unless a perturbation in the system is applied.

**Mathematical derivation of the molecular incompressibility**

We will assume that there is no void space or discontinuities in the gel. Therefore, the differential volume of the gel is given by equation (2).

$$dV_{gel} = dV_{network} + dV_{solvent} \tag{2}$$

Expressing the volume of the gel and the network with respect to the reference state (Lagrangian point of view) and noticing that $V^0$ refers to the volume at time $t=0$ and $V$ refers to the volume at time $t \geq 0$, we get

$$det(\boldsymbol{F}) dV^0_{gel} = dV^0_{network} + dV_{gel} \tag{3}$$

Also, note that $dV^0_{network} = dV^0_{gel}$ (see, i.e., figure 3) and we have used the incompresibility of the dry network. Dividing equation (3) by $dV^0_{network}$ and noticing that $dV_{solvent}/dV^0_{network}$ is the volume of the small molecules (solvent) in the gel divided by the volume of the dry network. We will approximate the latter, using $vC$. Where $v$ is the volume of the solvent molecules and $C$ is a scalar field that describes the distribution of the solvent molecules in the gel. With these assumptions, we get the molecular incompressibility condition given in equation (1).

**3. FLORY-REHNER FREE ENERGY FUNCTION**

For the purpose of this analysis we will work with the Flory-Rehner free energy function [9], defined by equation (4)

$$w(\boldsymbol{F},C) = \frac{1}{2} NkT [\boldsymbol{F}:\boldsymbol{F} - 3 - 2\log(det\,\boldsymbol{F})] - \frac{kt}{v}\left[vC\log(1+\frac{1}{vC}) + \frac{\chi}{1+vC}\right] \quad (4)$$

Where $N$ is the number of polymeric chains per reference volume, $kT$ is the absolute temperature in units of energy, $\chi$ is a dimensionless measure of the enthalpy of mixing. In equation (4) we have two independent variables, $\boldsymbol{F}$ and $C$. Therefore, if we want to use equation (4) in a variational principle, then we will need to use a two-field variational principle [13]. Other way is to use a Lagrange transform and the molecular incompressibility to transform equation (4) into a function of the deformation gradient only. When the gel equilibrates with the solvent and the mechanical load, the *chemical potential of the solvent molecules*, μ, is homogeneous in the external solvent and inside the gel, as is presented in figure 4. In thermodynamics, the chemical potential is defined by equation (5), keeping entropy, volume and number of particles constant.

$$\mu = \frac{\partial W(\boldsymbol{F},C)}{\partial C} \quad (5)$$

Using equation (5), the Legendre transform and the molecular incompressibilty, we get equation (6). Later we will replace $C$ using the molecular incompressibility condition.

$$\hat{W}(\boldsymbol{F},\mu) = W(\boldsymbol{F},C) - \mu C \quad (6)$$

A combination of equation (4), equation (6) and equation (1) gives the desired free energy function.

$$\hat{W}(\boldsymbol{F},\mu)=\frac{1}{2}NkT(I-3-2\log J)-\frac{kT}{v}\left[(J-1)\log(\frac{J}{J-1})+\frac{\chi}{J}\right]-\frac{\mu}{v}(J-1) \qquad (7)$$

Where $I=\boldsymbol{F}:\boldsymbol{F}$ is the first invariant of the right Cauchy-Green deformation tensor. Seems like there is an inconsistency in equation (7) when $J=1$, that is, solvent free, but it is not true because equation (7) is bounded, in other words, $\lim_{F=I}\hat{W}(\boldsymbol{F},\mu)=0$.

## 4. HYPERELASTIC CONSTITUTIVE EQUATION FOR GELS

We will consider the hydrogel as a hyperelastic material, therefore the constitutive equation is given by equation (8), where $\boldsymbol{P}$ is known as the first Piola-Kirchhoff stress tensor [13] or nominal stress.

$$\boldsymbol{P}=\frac{\partial \hat{W}(\boldsymbol{F},\mu)}{\partial \boldsymbol{F}} \qquad (8)$$

Next step is to replace equation (7) into equation (8) and take the derivative. In the following, we are going to use equation (9)

$$\frac{\partial J}{\partial \boldsymbol{F}}=J\boldsymbol{F}^{-T} \qquad (9)$$

For simplicity we will break the derivative of the free energy function in three parts.

Part 1

$$\frac{\partial}{\partial \boldsymbol{F}}\left[\frac{1}{2}NkT(\boldsymbol{F}:\boldsymbol{F}-3-2\log J)\right]=NkT(\boldsymbol{F}-\boldsymbol{F}^{-T})$$

Part 2

$$\frac{kT}{v}\frac{\partial}{\partial \boldsymbol{F}}\left[(J-1)\log(\frac{J}{J-1})+\frac{\chi}{J}\right]=\frac{kT}{v}\left[J\log(\frac{J}{J-1})-1-\frac{\chi}{J}\right]\boldsymbol{F}^{-T}$$

Part 3

$$\frac{\partial}{\partial \boldsymbol{F}}\left[\frac{\mu}{v}(J-1)\right]=\frac{\mu}{v}J\boldsymbol{F}^{-T}$$

Replacing part 1, part 2 and part 3 into equation (8), we get the constitutive equation for this particular problem.

$$\boldsymbol{P}=NkT(\boldsymbol{F}-\boldsymbol{F}^{-T})-\frac{kT}{v}\left[J\log(\frac{J}{J-1}-1-\frac{\chi}{J})\right]\boldsymbol{F}^{-T}-\frac{\mu}{v}J\boldsymbol{F}^{-T} \qquad (10)$$

Note that the first Piola-Kirchhoff stress tensor is singular when the network is solvent free. Therefore, we can not use equation (10) in the reference state the dry network. In the following section is proposed to use a free-swelling as a reference state.

**Free swelling as a reference state**

To avoid the singularity in the free-energy function, we choose as our reference state the free swelling state [9]. We will need to specify equation (10) in terms of this state. To do it, we will use a multiplicative product of the deformation gradient, as in figure 5.

With reference to figure 5 and with the definition of the deformation gradient we get

$$d\boldsymbol{x}'=\boldsymbol{F}_0 d\boldsymbol{X}; \quad d\boldsymbol{x}=\boldsymbol{F}d\boldsymbol{x}'; \quad d\boldsymbol{x}=\boldsymbol{F}'\boldsymbol{F}_0 d\boldsymbol{X}$$

where in the last step we have replaced $d\boldsymbol{x}'$ by $\boldsymbol{F}_0 d\boldsymbol{X}$ and $d\boldsymbol{x}'$ is a differential vector in the free swelling state, $d\boldsymbol{x}$ is a differential vector in the swelling state and $d\boldsymbol{X}$ is a differential vector in the reference state (dry network). Therefore, the multiplicative product of the deformation gradient is given by equation (11)

$$F = F' F_0 \tag{11}$$

The motion between the dry network state and the free swelling state is given by equation (12)

$$x' = \lambda_0 X \tag{12}$$

where $\lambda_0$ is the free swelling stretch. Using equation (12) and the definition of deformation gradient, we get the deformation gradient between the free swelling state and the dry network state, equation (13), where $I$ is the identity matrix.

$$F_0 = \frac{\partial x'}{\partial X} = \lambda_0 I \tag{13}$$

**Free swelling stretch**

For the evaluation of the free swelling stretch, we will use equation (10) and the fact that in the free swelling state $P = 0$, where $0$ is the zero matrix and $F = F_0$.

$$F_0^{-T} = \frac{1}{\lambda_0} I; \quad J = \lambda_0^3$$

with the above in mind, we get

$$0 = NkT(\lambda_0 - \frac{1}{\lambda_0})I - \frac{kT}{v}\left[\lambda_0^3 \log(\frac{\lambda_0^3}{\lambda_0^3 - 1}) - 1 - \frac{\chi}{\lambda_0^3}\right](\frac{1}{\lambda_0})I - \frac{\mu_0}{v}\lambda_0^3(\frac{1}{\lambda_0})I$$

It reduce to one equation with one unknowns

$$Nv\left(\frac{1}{\lambda_0}-\frac{1}{\lambda_0^3}\right)+\log\left(1-\frac{1}{\lambda_0^3}\right)+\frac{1}{\lambda_0^3}+\frac{\chi}{\lambda_0^6}=\frac{\mu_0}{kT} \tag{14}$$

Equation (14) is a nonlinear equation and could be easily solve using a Newton-Rahpson algorithm. Note that we have used the property: $\log(u/v)=-\log(v/u)$.

**Free energy as a function of the free swelling state**

We desire to reformulate equation (7) as a function of the free swelling state, that is, $F'$ and $\mu$. Before proceeds, it is important to say that the free energy function is defined per *unit reference volume*. To developed the relation we will use a dimensional analysis. In the following, $U$ is the strain energy in the swelling state, $V_0$ is the reference volume (dry network) and $V'$ is the volume in the free swelling state. Then,

$$\hat{W}=\frac{U}{V_0}; \quad \hat{W}'=\frac{U}{V'}$$

Note that the strain energy should be the same for both cases, the difference in the free energy function is due to the choice of the reference state. But, we can relate $V_0$ with $V'$ by mean of the *Jacobian*. Also, the Jacobian for the free swelling state is constant.

$$dV'=J_0 dV_0; \quad \int_{V'} dV'=J_0\int_{V_0} dV_0; \quad V'=J_0 V_0$$

Noticing that $J_0=\lambda_0^3$, then,

$$\hat{W}'(F',\mu)=\lambda_0^{-3}\hat{W}(F,\mu) \tag{15}$$

*Free energy*

The free energy as a function of the free swelling state is obtained using equation (7), equation

(15) and the expression obtained in previous sections.

$$\hat{W}'(\boldsymbol{F}',\mu) = \frac{\lambda_0^{-3}}{2} NkT \left[ \lambda_0^3 I' - 3 - 2\log(\lambda_0^3 J') \right] - \frac{kT}{v} \left[ (J' - \lambda_0^{-3}) \log\left(\frac{\lambda_0^3 J'}{\lambda_0^3 J' - 1}\right) + \frac{\chi}{\lambda_0^6 J'} \right] \\ \frac{-\mu}{v}(J' - \lambda_0^{-3}) \tag{16}$$

In the free swelling state, $J'=1$, but it does not cause problem, unless, $\lambda_0=1$, that is, the gel is incompressible, therefore, violating our assumption of compressibility.

**Constitutive equation in terms of invariant and elasticity tensor**

In this section is introduced the most general form of a stress relation and elasticity tensor in terms of the strain invariant for the Flory-Rehner free energy function as is applied to hyperelastic material models. The most general expression for the evaluation of stresses using equation (16) is

$$\boldsymbol{S} = 2\frac{W(\boldsymbol{C})}{\partial \boldsymbol{C}} = 2\left[\left(\frac{\partial W}{\partial I_1}\right)\boldsymbol{I} + I_3 \frac{\partial W}{\partial I_3}\boldsymbol{C}^{-1}\right] \tag{17}$$

where $\boldsymbol{S}$ is the second Piola-Kirchhoff stress tensor $I_1$ and $I_3$ are the invariants of $\boldsymbol{C}$ and $\boldsymbol{I}$ is the identity matrix. For the description of isotropic hyperelastic materials at finite strain considering equation (16) in the coupled form, with invariant, $I_1$ and, the $I_3$ most general elasticity tensor is

$$\boldsymbol{D} = 2\frac{\partial \boldsymbol{S}(\boldsymbol{C})}{\partial \boldsymbol{C}} = 4\frac{\partial^2 \hat{W}(I_1, I_3)}{\partial \boldsymbol{C} \partial \boldsymbol{C}} \tag{18}$$

$$D = \delta_1 C^{-1} \otimes C^{-1} + \delta_2 C^{-1} \odot C^{-1} \tag{19}$$

where

$$\delta_1 = 4(I_3 \frac{\partial \hat{W}'}{\partial I_3} + I_3^2 \frac{\partial^2 \hat{W}'}{\partial I_3 \partial I_3}) \tag{20}$$

$$\delta_2 = -4 I_3 \frac{\partial \hat{W}'}{\partial I_3} \tag{21}$$

$$(C^{-1} \otimes C^{-1})_{ijkl} = C^{-1}_{ij} C^{-1}_{kl} \tag{22}$$

$$(C^{-1} \odot C^{-1})_{ijkl} = \frac{1}{2}(C^{-1}_{ik} C^{-1}_{jl} + C^{-1}_{il} C^{-1}_{jk}) \tag{23}$$

*Notes about derivatives of the Flory-Rehner free energy function respect to invariants of Green Lagrange strain tensor.*

All the equation on this section assumed that the Flory-Rehner free-energy function is a function of the first and thirds invariants of $C$. However, equation (16) is expressed as a function of $I_1$ and $J = det\, F = I_3^{1/2}$. Therefore, to differentiate $W$ respect to $I_3$, we should use chain rule.

$$\frac{\partial \hat{W}'}{\partial I_3} = \frac{\partial \hat{W}'}{\partial J} \frac{\partial J}{\partial I_3} \tag{24}$$

$$\frac{\partial J}{\partial I_3} = \frac{1}{2} I_3^{-1/2} = \frac{1}{2} J^{-1} \tag{25}$$

$$\frac{\partial^2 J}{\partial I_3 \partial I_3} = \frac{-1}{4} I_3^{-3/2} = \frac{-1}{4} J^{-3} \tag{26}$$

$$\frac{\partial^2 \hat{W}'}{\partial I_3 \partial I_3} = \frac{\partial^2 \hat{W}'}{\partial J^2}\left(\frac{\partial J}{\partial I_3}\right)^2 + \frac{\partial \hat{W}'}{\partial J}\frac{\partial^2 J}{\partial I_3 \partial I_3} \tag{27}$$

# 5. VARIATIONAL FORM OF THE EQUATION OF MOTION AND FINITE ELEMENT EQUATIONS

The weak form of the boundary-value problem is also known as the principle of virtual work in mechanics of materials. In this method we want to express the boundary-value problem in a variational form (integral form) instead of the original differential form. Therefore, the requirements of the continuity of the trial solution is lower than in the strong form, it suggest the name weak form. For a more detail explanation and a proof that both formulations are equivalent, see [15].

A formal statement of the weak form of the boundary-value goes as follows: Given $\boldsymbol{b}:\Omega\to\Re^n, \bar{\boldsymbol{u}}:\Gamma_u\to\Re^n$ and $\bar{\boldsymbol{t}}:\Gamma_t\to\Re^n$, find $\boldsymbol{u}\in V$ such that for all $\eta\in W$,

$$\int_\Omega \sigma:\nabla\eta\, d\Omega - \int_{\Gamma_t}\bar{\boldsymbol{t}}\cdot\eta\, d\Gamma - \int_\Omega \boldsymbol{b}\cdot\eta\, d\Omega = 0 \tag{28}$$

where $\sigma$ is a second order tensor defined in terms of $\boldsymbol{u}$ and $V$ is the space for all the test functions (solution space) and $W$ is the space for all the trial functions (variation), such that the following is true,

$$V = \left\{\boldsymbol{u} \mid \boldsymbol{u}\in C^0(\boldsymbol{X}),\ \boldsymbol{u}=\bar{\boldsymbol{u}}\ on\ \Gamma_u\right\} \tag{29}$$

$$w = \{\eta \mid \eta \in C^0(X), \quad \eta = 0 \quad on \quad \Gamma_u\} \tag{30}$$

Figure 6 is a pictorial representation of the trial function. Note that $\eta$ is arbitrary and the only condition is that must be zero when evaluated at the boundary where the displacement has been specified. We could look upon $\eta$ as the virtual displacement field $\delta u$, defined on the current configuration. Taking this into account we want to recast equation (28) as a function of the variation of the displacement and the Cauchy stress tensor ($\sigma$). Considering the symmetry of the Cauchy stress tensor and

$$\delta e = \frac{1}{2}((\nabla \delta u)^T + \nabla \delta u) = sym(\nabla \delta u) \tag{31}$$

Where $\delta e$ is the variation of the Euler-Almansi strain tensor [13] defined in the current state. Substituting equation (31) into equation (28) we get,

$$\int_\Omega \sigma : \delta e \, d\Omega = \int_{\Gamma_t} \bar{t} \cdot \delta u \, d\Gamma + \int_\Omega b \cdot \delta u \, d\Omega \tag{32}$$

Equation (32) represent the virtaul desplacement principle's and it clearly state that the internal virtual work should equilibrate the external virtual work. This equation is true for small and large deformation and small or large stress ans strain. The use of equation (32) require the specification of a proper material model. The expression of the virtual work in a nonlinear regime has a fundamental difficulty that in most cases does not allow to use it, as is stated in equation (32).

This is because, in general, we do not know the domain of integration at the current time and the material model used in this paper is nonlinear.

**Total Lagrangian formulation**

In this paper we used the Total Lagrangian (TL) formulation to deal with the unknown integration domain. A TL formulation is based on the assumption that there exist an inverse mapping between the current and a reference state, that allows us, to refers everything to a known reference state at the beginning of the motion. To express the principle of virtual work in the reference state we will use the second Piola-Kirchhoff stress tensor and the Green-Lagrange strain tensor ( $E$ ).

$$\int_{\Omega} \sigma : \delta e \, d\Omega = \int_{\Omega_0} S : \delta E \, \delta \Omega \qquad (33)$$

Equation (33) is the basic expression of the TL formulation used in this paper. Note that the integration is respect to the known reference configuration at time $t=0$. It is equivalent as what we do in finite element analysis with the use of a parent element.

**Integration in a finite element solution based on Total Lagrangian approach**

Integrals on the reference configuration are related to integrals over the parent element by

$$\int_{\Omega_0} g(X) \, d\Omega_0 = \int_{\xi} g(\boldsymbol{\xi}) J_\xi \, d\xi \qquad (34)$$

where the integration is carried out on the parent element domain ( $\xi$ ) and $J_\xi$ is the Jacobian between the reference configuration and the parent element, as shown in figure 7 for a two dimensional quadrilateral element. Note that, $J_\xi$ does not depends o the current configuration and therefore just need to be evaluated at the beginning of the calculations.

**Finite element approximation**

To solve approximately the equations on previous section, we will use an isoparametric finite element discretization. The problem domain is defined by a mesh. Geometry and independent variables are interpolated using compact piecewise linear functions. In general, the geometry and variables will be interpolated using equation (35) for each element in the mesh.

$$\begin{aligned} \boldsymbol{x}(\xi) &= \boldsymbol{N}(\xi)\hat{\boldsymbol{x}} \\ \boldsymbol{u}(\xi) &= \boldsymbol{N}(\xi)\hat{\boldsymbol{u}} \end{aligned} \quad (35)$$

where $\xi$ corresponds to the parent element coordinate, $\boldsymbol{N}(\xi)$ are the standard finite element shape functions [19] and $\hat{\boldsymbol{x}}$ and $\hat{\boldsymbol{u}}$ are the vectors of coordinates nodal values and nodal variables, respectively. Then, using the previous interpolations and Voigt notation [16] we can express equation (33) in matrix form and we will call it, the internal nodal force vector.

$$\boldsymbol{F}_I = \int_{\Omega_0} \boldsymbol{B}_{0I}^T \{\boldsymbol{S}\} d\Omega_0 \quad (36)$$

where the $B_{0I}$ express the variation of the Green strain [16, 17] and is given as follows for a three dimensional problem,

$$B_{0I} = \begin{bmatrix} \frac{\partial N_I}{\partial X}\frac{\partial x}{\partial X} & \frac{\partial N_I}{\partial X}\frac{\partial y}{\partial X} & \frac{\partial N_I}{\partial X}\frac{\partial z}{\partial X} \\ \frac{\partial N_I}{\partial Y}\frac{\partial x}{\partial Y} & \frac{\partial N_I}{\partial Y}\frac{\partial y}{\partial Y} & \frac{\partial N_I}{\partial Y}\frac{\partial z}{\partial Y} \\ \frac{\partial N_I}{\partial Y}\frac{\partial x}{\partial Z}+\frac{\partial N_I}{\partial Z}\frac{\partial x}{\partial Y} & \frac{\partial N_I}{\partial Y}\frac{\partial y}{\partial Z}+\frac{\partial N_I}{\partial Z}\frac{\partial y}{\partial Y} & \frac{\partial N_I}{\partial Y}\frac{\partial z}{\partial Z}+\frac{\partial N_I}{\partial Z}\frac{\partial z}{\partial Y} \\ \frac{\partial N_I}{\partial X}\frac{\partial x}{\partial Y}+\frac{\partial N_I}{\partial Y}\frac{\partial x}{\partial X} & \frac{\partial N_I}{\partial X}\frac{\partial y}{\partial Y}+\frac{\partial N_I}{\partial Y}\frac{\partial y}{\partial X} & \frac{\partial N_I}{\partial X}\frac{\partial z}{\partial Y}+\frac{\partial N_I}{\partial Y}\frac{\partial z}{\partial X} \end{bmatrix} \tag{37}$$

and the second Piola-Kirchhoff stress "vector" is given by,

$$\{S\} = \begin{bmatrix} S_{11} & S_{22} & S_{33} & S_{23} & S_{13} & S_{12} \end{bmatrix}^T \tag{38}$$

Applying a linearization process to equation (36) and using finite element interpolation, we arrived to the following system of equations,

$$(K_{mat} + K_{geo})\Delta\hat{U} = R - F \tag{39}$$

Where $K_{mat}$ and $K_{geo}$ are the well known material and geometric stiffness matrix, respectively and $R$ is the vector of external point forces. This matrices are defined by,

$$K_{IJ}^{mat} = \int_{\Omega_0} B_{0I}^T D B_{0J} \, d\Omega \tag{40}$$

$$K_{IJ}^{geo} = I\int_{\Omega_0} B_{0I}^T S B_{0J} \, d\Omega \tag{41}$$

where $B_{0I} = \partial N_I / \partial X$, $D$ is the algorithmic tangent modulus ( equation (19) ) in Voigt notation and $I$ is the identity matrix.

## 6. RESULTS

**Free swelling**

Free swelling refers to a gel configuration in which there are no constraints for the motion. To simulate it into a finite element solution we removed the rigid body motion of the body, applying proper boundary conditions. In figure 8 is shown a comparison between the analytical and finite element results. The gel was simulated using one element as shown in figure 9.

Note also that the displacement field shown in figure 9 was measured respect to the free swelling state, because of our implementation of the Flory-Rehner free energy function. According to that, the relation used to calculate the free swelling stretch is given by equation (42)

$$\lambda = (1 + \frac{\Delta}{L})\lambda_0 \qquad (42)$$

where $\lambda_0$ is the initial free swelling stretch, $\Delta$ is the displacement, measure from the initial swelling state, and $L$ is the original free swelling cube length.

**A bar of a gel in contact with a solvent and subjected to a uniform axial stress**

In this example we put in contact a bar of gel with a solvent and an axial force was applied. The longitudinal stretch is given by $\lambda_1$ and the traverse stretches are given by $\lambda_2 = \lambda_3$. The analytical solution is given in equation (43) and it is used to calculate (using Newton's method) $\lambda_2$ for every known value $\lambda_1$. The bar was simulated using one finite element as shown in figure 10. The comparison between the analytical and numerical solution is presented in figure 11, and it shows that our material model implementation is correct.

$$Nv\left(\lambda_2 - \frac{1}{\lambda_2}\right) + \left[\lambda_1 \lambda_2^2 \log\left(1 - \frac{1}{\lambda_1 \lambda_2^2}\right) + 1 + \frac{\chi}{\lambda_1 \lambda_2^2} - \frac{\mu}{kT} \lambda_1 \lambda_2^2\right] \frac{1}{\lambda_2} = 0 \qquad (43)$$

## 7. CONCLUSIONS

We have shown some preliminaries results that involves the use of a Flory-Rehner free energy function in a Total Lagrangian framework. The results have shown that our implementation of the material model is working as expected. Further analysis that includes non homogeneous field should be done. The extremely large deformation in the gel imposes a difficult convergence condition to the finite element solution.

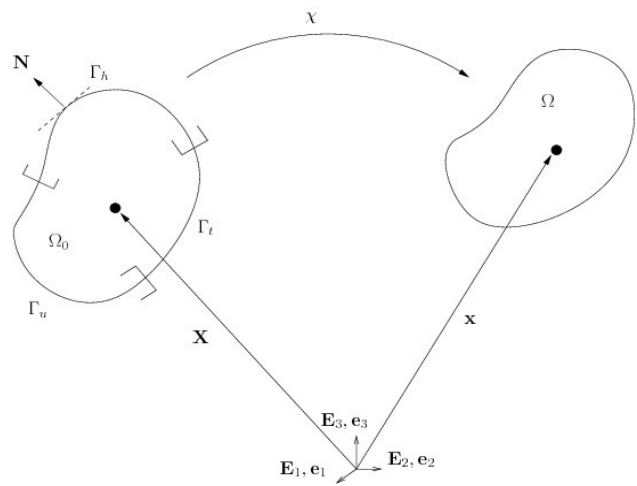

**Figure 1.** Reference and current configuration.

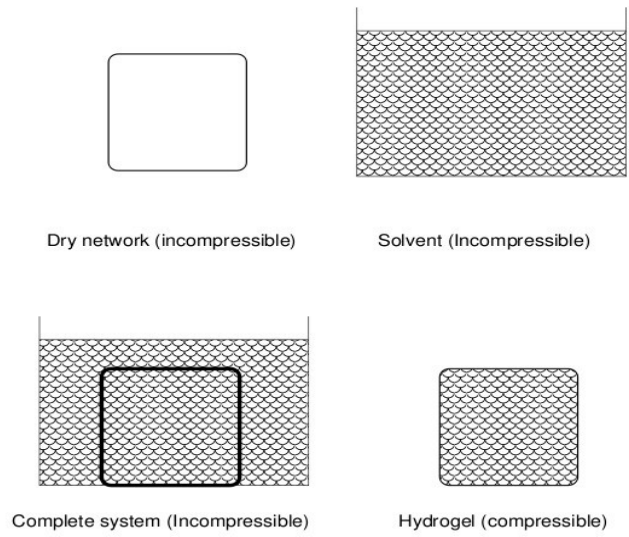

**Figure 2**. Pictorial definition of a hydrogel.

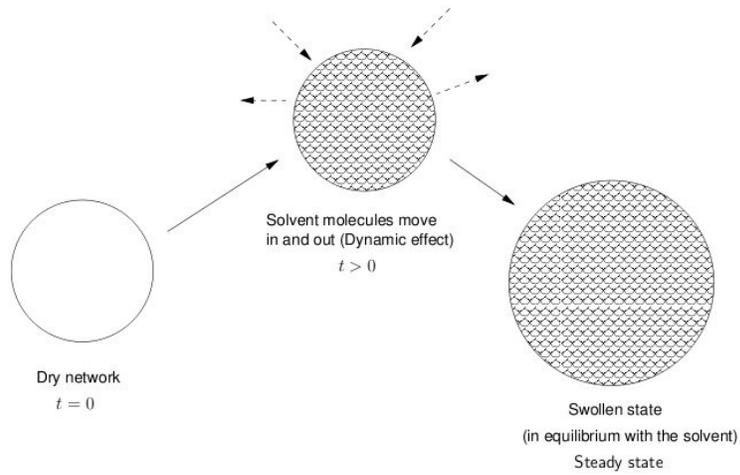

**Figure 3**. Pictorial of the kinetic of swelling gels.

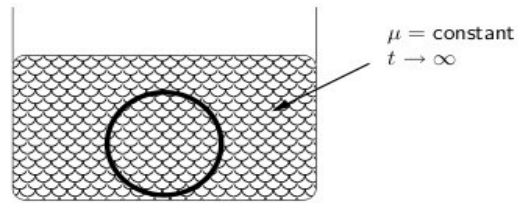

**Figure 4**. In equilibrium the chemical potential is constant for the solvent and the gel.

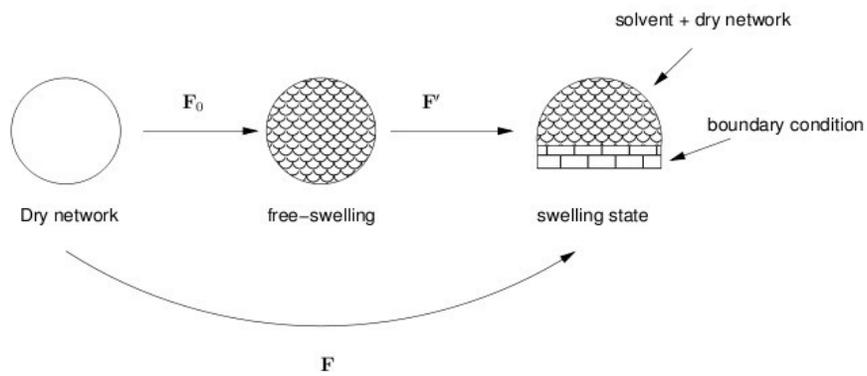

**Figure 5**. Multiplicative product of the deformation gradient.

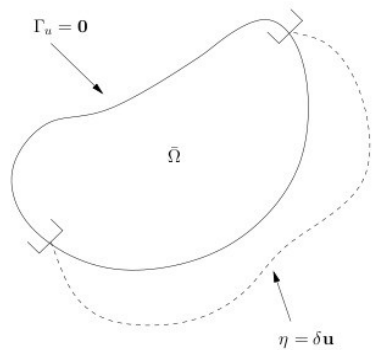

**Figure 6**. Pictorial definition of the trial function. The dash line corresponds to the trial function

respect to the equilibrium position.

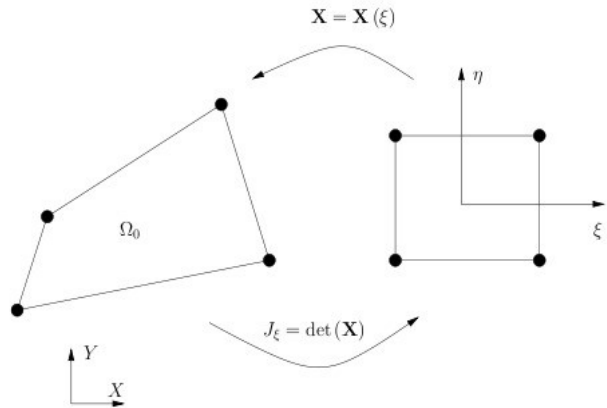

**Figure 7**. Mapping between the reference configuration and the parent element.

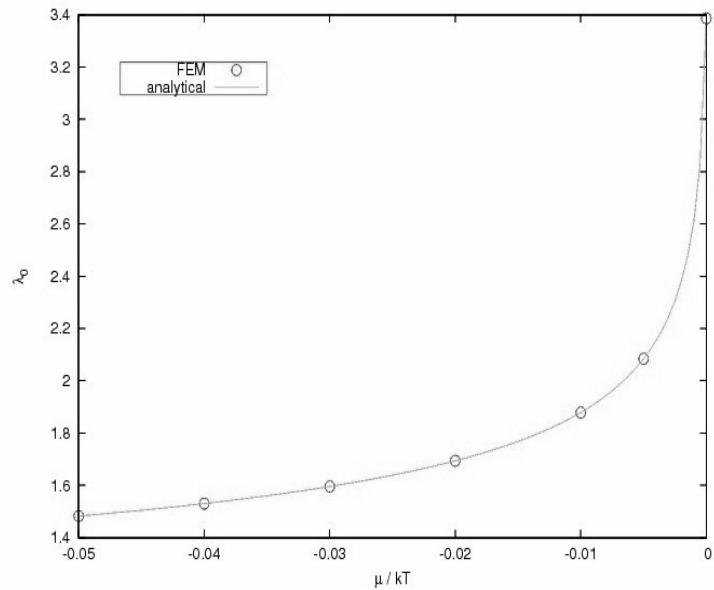

**Figure 8.** Free swelling response as a function of the normalized chemical potential of the solvent. *Parameters*:   $Nv=10^{-3}$   and   $\chi=0.1$   .

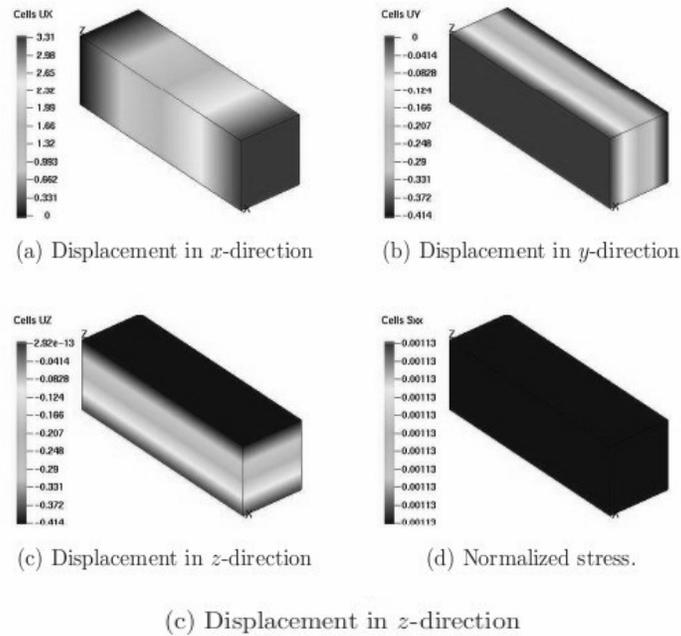

**Figure 9.** Displacement field for $\mu/kT=0$ and measure from a free swelling state at $\lambda_0=1.482$. *Parameters:* $Nv=10^{-3}$ and $X=0.1$ and cube with free swelling length $L=2$. **Figure 10.** Displacement field and uniaxial stress for $\mu/kT=0$ and measure from a free swelling state at $\lambda_0=3.390$. *Parameters:* $Nv=10^{-3}$ and $X=0.1$ and cube with free swelling length $L=2$.

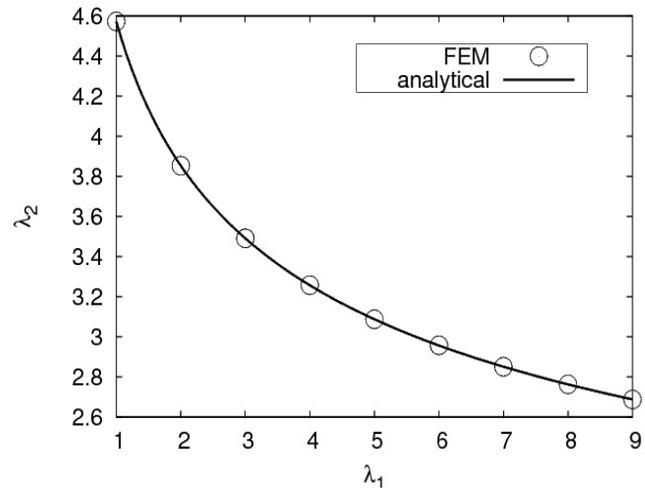

**Figure 11.** A bar of a gel is brought into contact with a solvent and subjected to a uniform axial stress. *Parameters*: $Nv=10^{-3}$ and $\chi=0.1$ .